\documentstyle[12pt]{article}
\begin{document}
\newcommand{\nc}{\newcommand}
\nc{\namelistlabel}[1]{\mbox{#1}\hfil}
\newenvironment{namelist}[1]{%
\begin{list}{}
{
\let\makelabel\namelistlabel
\settowidth{\labelwidth}{#1}
\setlength{\leftmargin}{1.1\labelwidth}
}
}{%
\end{list}}
\def\theequation{\thesection.\arabic{equation}}
\newtheorem{theorem}{\bf Theorem}[thetheorem]
\newtheorem{prop}{\bf Proposition}[theprop]
\newtheorem{corollary}{\bf Corollary}[thecorollary]
\newtheorem{remark}{\bf Remark}[theremark]
\newtheorem{lemma}{\bf Lemma}[thelemma]
\def\thetheorem{\thesection.\arabic{theorem}}
\def\theprop{\thesection.\arabic{prop}}
\def\theremark{\thesection.\arabic{remark}}
\def\thecorollary{\thesection.\arabic{corollary}}
\def\thelemma{\thesection.\arabic{lemma}}
\nc{\bsp}{\begin{sloppypar}}
\nc{\esp}{\end{sloppypar}}
\nc{\be}{\begin{equation}}
\nc{\ee}{\end{equation}}
\nc{\beanno}{\begin{eqnarray*}}
\nc{\inp}[2]{\left( {#1} ,\,{#2} \right)}
\nc{\dip}[2]{\left< {#1} ,\,{#2} \right>}
\nc{\disn}[1]{\|{#1}\|_h}
\nc{\pax}[1]{\frac{\partial{#1}}{\partial x}}
\nc{\tpar}[1]{\frac{\partial{#1}}{\partial t}}
\nc{\xpax}[2]{\frac{\partial^{#1}{#2}}{\partial x^{#1}}}
\nc{\pat}[2]{\frac{\partial^{#1}{#2}}{\partial t^{#1}}}
\nc{\ntpa}[2]{{\|\frac{\partial{#1}}{\partial t}\|}_{#2}}
\nc{\xpat}[2]{\frac{\partial^{#1}{#2}}{\partial t \partial x}}
\nc{\npat}[3]{{\|\frac{\partial^{#1}{#2}}{\partial t^{#1}}\|}_{#3}}
\nc{\xkpat}[3]{\frac{\partial^{#1}{#2}}{\partial t^{#3} \partial x}}
\nc{\jxpat}[3]{\frac{\partial^{#1}{#2}}{\partial t \partial x^{#3} }}
\nc{\eeanno}{\end{eqnarray*}}
\nc{\bea}{\begin{eqnarray}}
\nc{\eea}{\end{eqnarray}}
\nc{\ba}{\begin{array}}
\nc{\ea}{\end{array}}
\nc{\nno}{\nonumber}
\nc{\dou}{\partial}
\nc{\bc}{\begin{center}}
\nc{\ec}{\end{center}}
 \nc{\bb}{\mbox{\hspace{.25cm}}}
\nc{\bite}{\begin{itemize}}
\nc{\eite}{\end{itemize}}
\nc{\bth}{\begin{theorem}}
\nc{\eth}{\end{theorem}}
\nc{\bpr}{\begin{prop}}
\nc{\epr}{\end{prop}}
\nc{\blem}{\begin{lemma}}
\nc{\elem}{\end{lemma}}
\nc{\benu}{\begin{enumerate}}
\nc{\eenu}{\end{enumerate}}
\nc{\bcor}{\begin{corollary}}
\nc{\ecor}{\end{corollary}}
\nc{\Pf}{{\bf Proof. }}
\nc{\ddt}{\f{d}{dt}}
\nc{\ddr}{\f{d}{dr}}
\nc{\ddx}{\f{d}{dx}}
\nc{\DDx}{\f{d^{2}}{dx^{2}}}
\nc{\dodot}{\f{\dou}{\dou t}}
\nc{\dodori}{\f{\dou}{\dou r_i}}
\nc{\dodorj}{\f{\dou}{\dou r_j}}
\nc{\N}{I\!\!N}
\nc{\R}{I\!\!R}
\nc{\Hy}{I\!\!H}
\nc{\C}{I\!\!\!\!C}
\nc{\na}{\nabla}
\nc{\la}{\lambda}
\nc{\s}{\sinh}
\nc{\co}{\cosh}
\nc{\vl}{V_{\la}}
\nc{\ga}{\gamma}
\nc{\vg}{V_{\ga}}
\nc{\T}{\Theta}
\nc{\th}{\theta}
\nc{\f}{\frac}
\nc{\rw}{\rightarrow}
\nc{\lrw}{\Longrightarrow}
\nc{\om}{\omega}
\nc{\Om}{\Omega}
\nc{\al}{\alpha}
\nc{\ro}{\varrho}
\nc{\qed}{\phantom{a}{\hfill \rule{2.5mm}{2.5mm}}}
\nc{\D}{\Delta}
\nc{\TD}{\tilde{\D}}
\nc{\si}{\sigma}
\nc{\noi}{\noindent}
\nc{\bib}{\bibitem}

\title{\bf Almost Diameter Rigidity for Cayley Plane} 
\author{Akhil Ranjan }
\maketitle
\begin{abstract}
In this paper we give a generalisation of the Radius Rigidity theorem
of F.Wilhelm. This is done by showing that if a Riemannian submersion of 
$S^{15}$ with 7-dimensional fibres has at least one fibre which is a 
great sphere then all the fibres are so. Some weaker than known conditions
which force the existence of such a fibre are also discussed.
\end{abstract}
{\bf Keywords:} Curvature, Diameter, Riemannian submersions, Isoparametric
foliations.
\footnotetext{AMS Subject Classification: 53C35}
\footnotetext{MRR 07-97}
\section{\bf Introduction}
Gromll and Grove in their beautiful papers \cite{GG1} and \cite{GG2}
almost classified the Riemannian submersions of the round spheres and applied
it to strengthen M. Berger's rigidity theorem characterising the compact
rank one symmetric spaces. They termed it as the Diameter Rigidity
Theorem. The analysis done in \cite{GG1} showed that the hypotheses 
involved in the Diameter Rigidity Theorem gave rise to a Riemannian submersion
of the unit tangent sphere at some point onto its cut-locus. The
classification of these submersions then enables one to prove symmetry.
However, just one case, that of Cayley-plane remained unsolved
as the corresponding classification of Riemannian submersions of $S^{15}$ 
with 7-dimensional fibres remained unsettled. Given such a submersion,
F. Wilhelm has shown that the set of points in the base 
over which the fibres are great spheres, has strong geometric
properties (Main Lemma in \cite{W}) and in conjunction with the assumed 
lower bound on the radius
shows that all the fibres, in the Riemannian submersion which arises, as in
\cite{GG1} are great spheres. This is what precisely is missing if only
diametrer condition is assumed. 

In the present work it is shown that Wilhelm's Main Lemma can be
considerably improved and consequently the following {\bf Almost Diameter
Rigidity} can be proved.
\begin{theorem}   
If a nonspherical Riemannian manifold $M$ has sectional curvatures $\geq
1$ and has an equilateral triangle of sides $\pi/2$, then it must be a
symmetric space.
\end{theorem}
In this paper we will only consider the case where $M$ is a cohomology
$CaP^2$.\\
It is interesting to note that the above can also be compared with
{\it corollary II} of \cite{Du} where it is required that every
pair of points $\pi/2$ distance apart be completed into an equilateral
triangle.

It is indicated towards the end of the paper how one can further weaken a 
little even this equilateral triangle condition.

To prove the above theorem we prove the following very strong version
of Wilhelm's Main Lemma
\begin{theorem}
Given a Riemannian submersion of $S^{15}$ with 7-dimensional fibres,
the set of those points in the base manifold $B$ over which the fibres 
are great spheres, is either empty or all of $B$.
\end{theorem}

\section{\bf Notations and Preliminaries}
Let us recall some standard notations and results used in this
paper. Let $S^{n+k}$ be the sphere of constant sectional curvature,
the constant being assumed to be unity. Let it be fibred
by leaves of dimension $k$ so that the metric is $"bundle-like"$ on
the sphere or in other words we have a Riemannian submersion onto a
Riemannian manifold $B$. O'Neill \cite{O'N}
introduced tensors $T$ and $A$ for any metrically foliated manifold $M$,
and called them {\it second fundamental tensor} and {\it integrability tensor}
respectively.First of all one decomposes the tangent bundle into a 
direct sum ${\cal V}\oplus{\cal H}$ where ${\cal V}$ is the set of
vectors tangential to the leaves and ${\cal H}$ are those normal to
the leaves. They are called vertical and horizontal vectors
respectively. For any tangent vector $e$,  we write $e = e^v + e^h$
or sometimes ${\cal V}(e) + {\cal H}(e)$,
uniquely as a sum of vertical and horizontal vectors. For $U$ and
$V$  smooth vertical vector fields and $X$ and $Y$  horizontal,
we define
\be
A_XY = (\na_XY)^v,\; {\rm and}\; A_XV = (\na_XV)^h
\ee
\be
T_VU = (\na_VU)^h,\;{\rm and}\; T_VX = (\na_VX)^v
\ee
Here $\na$ denotes the Levi-Civita connection on the manifold.
For any horizontal vector $x$ and vertical vector $v$ at a point
$p\in M,\,A_x {\rm and}\,T_v$ are skewsymmetric endomorphisms of
$T_pM$ which interchange vertical and horizontal vectors. Moreover
we have
\be
A_xy = -A_yx,\; {\rm and}\; T_uv = T_vu
\ee
For convenience we also write
\be
A^vx = A_xv,\;{\rm and}\; T^xv = T_vx
\ee
This makes $A^v$ (for a vertical $v$) act as a skewsymmetric operator
on horizontal vectors while $T^x$ (for $x$ horizontal) becomes a
symmetric endomorphism of vertical vectors. The latter is written 
as $S_x$ in \cite{GG2}.It is the second fundamental form operator.

Finally, we make $\na^v$ denote the (Riemannian) connection on ${\cal V}$
obtained by vertical projection of $\na$ and $\na^h$ that on ${\cal H}$.
Also as in \cite{O'N}, a horizontal vector field $X$ which projects to a
well defined vector field in $B$, will be referred to as {\it basic}.
It satisfies the differential equation
\be
\na^h_VX = A_XV = A^VX
\ee
where as usual $V$ is any vertical field.
\section{\bf Some Computations}
\begin{lemma}
Let $X$ and $Y$ be basic fields along a fibre of a Riemannian submersion
then
\be
div(A_XY) = -\sum_i<(\na_{v_i}A)^{v_i}X,Y>
\ee
In the above equation, $div$ denotes the divergence operator in the 
{\it fibre} and accordingly $\{v_i\}$ is an orthonormal basis of 
vertical vectors at any point of the fibre where divergence is to be found.
\end{lemma}
{\bf Proof:} Straightforward. The terms involving derivatives of $X$
and $Y$ cancel out on invoking skew-symmetries of $A$.\qed
\begin{corollary}
In the case of a Riemannian submersion of the round sphere, the fields
$A_XY$ are divergence free for $X$ and $Y$ basic along a fibre.
\end{corollary}
{\bf Proof:} From equation $\{2\}$ of \cite{O'N} and using isoparametricity
\cite{GG2} we easily see that $div(A_XY)$ is constant along the fibre.
Since average value of divergence of a tangent vector field on a closed
manifold is zero, we get the result.\qed

{\bf Remark:} Isoparametricity shows that the terms $<T^X,T^Y>$ and
$<\na_XK,Y>$ are constant along the fibre, where $K$ is the mean
curvature vector field. For a simple geometric proof of isoparametricity 
one can also see \cite{AR}. Also the constant sectional curvature implies
that $|A_XY|^2$ is constant along the fibres for basic fields $X$ and 
$Y$ (see \cite{GG2} or use eqn.\{4\} of \cite{O'N}).
\begin{corollary}
Let $X$ be a basic field along a fibre of a Riemannian submersion of
the round sphere, then
\be
\Delta X = \sum_i(A^{v_i})^2X.
\ee
Here $\Delta$ is the Laplacian computed by using the (Riemannian)
connection $\na^h$ on ${\cal H}$ restricted to the fibre.
\end{corollary}
{\bf Proof:} 
\be
\Delta X = \sum_i[(A^{v_i})^2X +(\na_{v_i}A)^{v_i}(X)]^h
         = \sum_i(A^{v_i})^2X
\ee
\hfill{\qed}
Now in the case of spheres the operator $\sum_i(A^{v_i})^2$ actually
descends to the base space $B$ and defines a symmetric negative definite
endomorphism of the tangent space $T_bB$ at each point $b\in B$. This is
due to last statement in the remark above.
If we fix one such $b$ over which our fibre under consideration lies,
we see that one can get an orthonormal framing of the fibre by basic
fields $\{X_1,X_2,...,X_n\}$ which are eigen vectors of the operator
$\sum_i(A^{v_i})^2$ and hence also of $\Delta$.

\section{Proofs of the Theorems}
We will be restricting ourselves to the case of a Riemannian submersion
of $S^{2k+1}$ onto $B^{k+1}$ with $k$ dimensional fibres one of which,
say $F$, is a great sphere. In this case its orthogonal complement
$F^{\perp}$ is also a great sphere and a fibre. We have a parallel
orthonormal framing $\{E_0,E_1,...,E_k\}$ of the horizontal vectors
along $F$. Let $\{X_j,0\leq j\leq k\}$ be the basic framing by the eigen 
vectors of the Laplacian and let, for simplicity of notation, $X$ be any 
one of these with eigen value $-\lambda$.
Clearly, one can write $X = \sum_j f_jE_j$ for suitable uniqely determined
smooth functions on $F$. This easily yields that each $f_j$ is an eigen
function of the Laplacian with the same eigen value $-\lambda$. Thus each
$f_i$ is a spherical harmonic in $k+1$ variables and $\lambda$ is in the
spectrum of $(S^k,g_{can})$.

Now fix this basic eigen vector field $X$ along $F$.Let $\{V_i,1\leq i\leq
k\}$ be an orthonormal framing of $TF$ by eigen vectors of the operator
$A_X^2$ acting on vertical vectors. This can be done as follows. As $T^X$
vanishes, $A_X$ is injective on ${\cal V}$ (see \cite{GG2},\cite{AR}) and 
therefore
nonsingular on $X^{\perp}$ due to dimension restrictions. Hence can find 
a basic $k$-frame $\{Y_i\}$ such that $A_X^2Y_i = -\mu_i^2Y_i$ 
with $\mu_i>0$. Set $V_i = A_XY_i/\mu_i$. Trivially, 
$A_X^2V_i = -\mu_i^2V_i$.     
\begin{lemma} As an operator on the vertical space at any point $p\in F$,
trace of $A_X^2$ satisfies the inequality
$$ trA_X^2 \leq -k $$.
\end{lemma}
{\bf Proof:} Let $\gamma$ be a horizontal geodesic starting at $p$ in
the direction $X(p)$. Let $v_i$ denote the parallel translate of $V_i(p)$
along $\gamma$ and $e_i$ that of $E_i = A_XV_i(p)$. Then the horizontal
holonomy displacement of $V_i(p)$ along $\gamma$ is given by
$$V_i(t) = \cos tv_i(t) + \sin te_i(t)$$
In particular, $V_i(\pi/2)$ is a basis of $T_{\gamma(\pi/2)}F^{\perp}$.
Consider the matrix $$M(t) = ((<V_i(t),V_j(t)>)) =
\cos^2tI+\sin^2t((\mu_i^2\delta_{ij})) = diag(\cos^2t+\mu_i^2\sin^2t) $$.
Due to isoparametricity, $\sqrt detM(t)$ gives the ratio of the volume
of the fibre $F(t)$ through $\gamma(t)$ and that of initial fibre $F$.
Since $F$ is isometric to $F(\pi/2) = F^{\perp}$ we get
$$ \prod_i\mu_i = 1$$
This implies $\sum_i\mu_i^2 \geq k$ as claimed. \qed
\begin{lemma} Let $X = \sum_if_iE_i$ for suitable harmonic polynomials
$\{f_i\}$ in $k+1$ variables $\{u_i,0\leq i\leq k\}$ and parallel 
orthonormal frame $\{E_i\}$ for any basic eigen vector field $X$ along
$F$, then each $f_i$ is linear.
\end{lemma}
{\bf Proof:} $\{u_i\}$ are coordinates in the $\R^{k+1}$ spanned by $F$.
Likewise, let $\{v_i\}$ denote the remaining coordinates in the Euclidean
space of $F^{\perp}$ The horizontal holonomy map from $F$ to $F^{\perp}$
induced by geodesics stating from various points of $F$ in the direction
of $X$ is given by
$$ v_i = f_i(u_0,...,u_k), 0\leq i\leq k $$
Thus it is a diffeomorphism which is algebraic. We will see presently that
its inverse is also algebraic.
The holonomy displacement produces a basic field $Y$ along $F^{\perp}$
and by the same considerations as for $F$, it is a restriction of a
polynomial vectorfield i.e. the coefficient functions with respect to
a parallel orthonormal framing of ${\cal H}$ are polynomials 
$g_i(v_0,...,v_k)$. But the holonomy displacement from $F^{\perp}$ to
$F$ via $Y$ is same as that given by $X$ at $t=\pi$. Hence
$g_i(f_0,...,f_k) = -u_i$. This proves that $f_0,...,f_k$ give an
algebraic equivalence. It follows therefore that $\sum_if_i^2 -1$
should generate the ideal $(\sum_iu_i^2 -1)$ and this forces $f_i$
to be linear homogeneous. \qed
\begin{corollary} As an operator on vertical vectors,
$$A_X^2 = -I$$ 
for any unit horizontal vector $X$ at any point of $F$.
\end{corollary}
{\bf Proof:} Since the coefficients of any basic field along $F$ are
linear it follows that $\lambda_i = k$ for every $i$.
Therefore, $\sum_i(A^{v_i})^2 = -kI$ and consequently for any unit
horizontal vector at any point of $F$, $\sum_i|A_Xv_i|^2 = k$.
This is same as $trA_X^2 = -k$ and this in turn forces $\mu_i = 1$ for 
each $i$. Clearly then $A_X^2 = -I$. \qed

{\bf Proof of Theorem 1.2:} From the above corollary we see that
$\sqrt detM(t) = 1$ for all $t$. Hence volume of $F(t)$ is constant
and equal to that of $F(0) = F$. Thus every fibre is of same volume 
as $vol(F)$ and therefore again by isoparametricity, each fibre is
a great $k$-sphere. Now it follows from the results proved in \cite{AR'}
that the submersion is congruent to the Hopf fibration. \qed

{\bf Proof of Theorem 1.1:} Let $x$, $y$, and $z$ be mutually distance
$\pi/2$ apart. In this situation $y$ and $z$ both are in the dual set
$\{x\}'$ of $x$. (See \cite{GG1},\cite{W} for details about dual sets.)
From \cite{GG1} we know that there is a Riemannian submersion
$$ \exp_x :S_x \rw \{x\}'$$
with 7-dimensional fibres from the unit tangent sphere at $x$, 
and similarly for $y$ and $z$. As argued 
in \cite{W} this forces at least one fibre to be totally geodesic
in each case. But then $\exp_x$ is congruent to the Hopf fibration
and this makes the space isometric to $(CaP^2, g_{can})$ as remarked
in \cite{GG1}.
\section{Concluding Remarks}
Due to isoparametricity of Riemannian submersions of spheres, it is easy
to see that the procedure of averageing the smooth functions over the 
fibres projects any eigen-space of the Laplacian orthogonally into
itself. Now the question which arises is whether the spherical harmonics
of degree two admit a nonzero function which is constant along the fibres
or not. If yes then the critical sets being great spheres, we get smaller
saturated great sheres. In case of $S^{15}$, if there is a fibre trapped
in the zone $\R^{15}\times [-1/4,1/4]$, then one can show that the harmonic
$x^2-1/16$ can be averaged to give a nonzero harmonic, $x$ being the last
coordinate in the above decomposition. Hence there exists  a proper
great sphere which is a union of fibres. Due to topological restrictions,
in case we are dealing with 7-dimensional fibres, it can only be a single
fibre. Thus the Main Theorem is further strengthened to
\begin{theorem} If a Riemannian submersion of $S^{15}$ is given with
7-dimensional fibres, and there is a fibre trapped in the zone
$\R^{15}\times [-1/4,1/4]$ then it must be congruent to a Hopf fibration.
\end{theorem}
\begin{corollary}
If $S^{15}\rw B^8$ is a Riemannian submersion such that the diameter
of the base $B$ is at least $\cos^{-1}(1/4)$, then it must be congruent
to the Hopf fibration.
\end{corollary} 
Now Theorem 1.1 can be rephrased suitably with a weaker hypothesis.

\noi
Department of Mathematics\\
Indian Institute of Technology\\
Mumbai 400076, INDIA\\
email: aranjan@ganit.math.iitb.ernet.in 
\end{document}